\documentclass[twocolumn]{autart}

\usepackage{amsfonts}
\usepackage{amssymb}
\usepackage{amsmath}
\usepackage{cite}
\usepackage{url}
\usepackage{graphicx}
\usepackage{subfigure}

\newtheorem{theorem}{Theorem}
\newtheorem{corollary}{Corollary}

\newtheorem{remark}{Remark}

\newcommand{\R}{\mathbb{R}}

\newcommand{\Rnm}[2]{\R^{{#1}\times{#2}}}
\newcommand{\bb}{(\cdot)}
\newcommand{\support}[2]{\operatorname{h}(#1,#2)}
\newcommand{\supportbb}[1]{\operatorname{h}(#1,\cdot)}
\newcommand{\gauge}[2]{\operatorname{g}(#1,#2)}
\newcommand{\gaugebb}[1]{\operatorname{g}(#1,\cdot)}

\newcommand{\mc}[1]{\mathcal{#1}}

\numberwithin{equation}{section}
\date{\today}
\begin{document}
\begin{frontmatter}
\runtitle{MLF: Alternative Characterization and Implicit Representation}
  \title{Minkowski--Lyapunov Functions: Alternative Characterization and Implicit Representation}
  \author[bit]{Sa\v{s}a~V.~Rakovi\'{c}\thanksref{corr}} 
	\address[bit]{Beijing Institute of Technology, Beijing, China}
  \thanks[corr]{E--mail: sasa.v.rakovic@gmail.com. Tel.: +447799775366.}
\begin{abstract}
An alternative characterization of  Minkowski--Lyapunov functions is derived. The derived characterization enables a computationally efficient utilization of Minkowski--Lyapunov functions in arbitrary finite dimensions. Due to intrinsic duality, the developed results apply in a direct manner to the characterization and utilization of robust positively invariant sets.
\end{abstract}
\end{frontmatter}

\section{Background}
\label{sec:01}
The characterization and computation of minimal and maximal robust positively invariant sets as well as their approximations are important research themes~\cite{kuntsevich:pshenichnyi:1996,kolmanovsky:gilbert:1998,rakovic:kerrigan:kouramas:mayne:2004b,limpiyamitr:ohta:2005,rakovic:2007,artstein:rakovic:2008}. A beneficial one to one correspondence between  robust positively invariant sets and Minkowski--Lyapunov functions has been established in a recent contribution~\cite{rakovic:2020.a}. The theoretical relevance of robust positively invariant sets and Minkowski--Lyapunov functions is amplified by their multifaceted practical utility. \emph{Inter alia}, robust positively invariant sets can be used as the target sets for robust time optimal controllers~\cite{mayne:schroeder:1997},  the uncertainty bounding sets for fast reference governors~\cite{gilbert:kolmanovsky:1999} and the tube cross--section shape sets for rigid tube model predictive controllers~\cite{mayne:seron:rakovic:2005}. Likewise, Minkowski--Lyapunov functions and their sublevel sets can be used as the terminal cost functions and constraint sets for consistently improving optimal control and stabilizing model predictive control~\cite{rawlings:mayne:2009}.  A more detailed insight into the theory, computation and applications of robust positively invariant sets and Minkowski--Lyapunov functions can be found in~\cite{kuntsevich:pshenichnyi:1996,kolmanovsky:gilbert:1998,rakovic:kerrigan:kouramas:mayne:2004b,limpiyamitr:ohta:2005,rakovic:2007,artstein:rakovic:2008,rakovic:2020.a,mayne:schroeder:1997,gilbert:kolmanovsky:1999,mayne:seron:rakovic:2005,blanchini:1999,blanchini:miani:2008,rawlings:mayne:2009,giesl:hafstein:2015} and numerous references therein. 

This note addresses an apparent lack of computational methods enabling a numerically efficient utilization of robust positively invariant sets and Minkowski--Lyapunov functions in arbitrary finite dimensions.  In Section~\ref{sec:02}, we derive alternative characterizations of Minkowski--Lyapunov functions and the fundamental Minkowski--Lyapunov function. In Section~\ref{sec:03}, we make use of these novel characterizations, in conjunction with implicit representations of  Minkowski functions, in order to create a potent platform for a computationally efficient and dynamically compatible utilization of Minkowski--Lyapunov functions in arbitrary finite dimensions; we also discuss a method for an alternative computation of the fundamental Minkowski--Lyapunov function. In Section~\ref{sec:04}, we show that, in light of intrinsic duality~\cite{rakovic:2020.a}, the derived results apply in a direct manner  to robust positively invariant sets.  In Section~\ref{sec:05}, we provide closing remarks including our numerical experience.  

\textbf{Basic Nomenclature.}
The spectral radius $\rho(M)$ of a matrix $M\in \Rnm{n}{n}$ is the largest absolute value of its eigenvalues. A matrix $M\in \Rnm{n}{n}$ is strictly stable if and only if $\rho(M)<1$.
 The Minkowski sum of nonempty sets $\mc{X}$ and $\mc{Y}$ in $\R^n$ is denoted by 
\begin{equation*}
\mc{X}\oplus \mc{Y}:=\{ x+y\ :\  x\in \mc{X},\ y\in \mc{Y}\}.
\end{equation*}
The image $M\mc{X}$ and the preimage $M^{-1}\mc{X}$ of a nonempty set $\mc{X}$ under a matrix of compatible dimensions (or a scalar) $M$ are denoted, respectively, by
\begin{equation*}
M\mc{X}:=\{Mx\ :\ x\in \mc{X}\}\text{ and }M^{-1}\mc{X}:=\{x\ :\ Mx\in \mc{X}\}.
\end{equation*} 
A $D$--set in $\R^n$ is a closed convex subset of $\R^n$ that contains the origin. A $C$--set in $\R^n$ is a bounded $D$--set in $\R^n$. A proper $D$--set in $\R^n$ is a closed convex subset of $\R^n$ that contains the origin in its interior. A proper $C$--set in $\R^n$ is a bounded proper $D$--set in $\R^n$. The Minkowski function $\gaugebb{\mc{X}}$ of a $D$--set $\mc{X}$ is given, for all $y\in\R^n$, by
\begin{equation*}
\gauge{\mc{X}}{y}:= \inf_\gamma \{\gamma\ :\ y\in \gamma \mc{X},\ \gamma\ge 0\}.
\end{equation*}

\textbf{Minkowski--Lyapunov Functions.} 
 A Minkowski--Lyapunov function~\cite{rakovic:2020.a} is the Minkowski function  $\gaugebb{\mc{S}}$ of a proper $C$--set $\mc{S}$ in $\R^n$ that verifies the Minkowski--Lyapunov inequality
\begin{equation*}
\forall x\in\R^n,\quad \gauge{\mathcal{S}}{Ax}+\gauge{\mathcal{Q}}{x}\le \gauge{\mathcal{S}}{x},
\end{equation*}
 in which $\mc{Q}$ is a given proper $C$--set in $\R^n$, and which is associated with the linear dynamics
\begin{equation*}
x^+=Ax,
\end{equation*}
where $x\in \R^n$ and $x^+\in \R^n$ are the current and successor states, and $A\in \Rnm{n}{n}$ is the state transition matrix.  The fundamental Minkowski--Lyapunov function~\cite{rakovic:2020.a} is the Minkowski function $\gaugebb{\mc{S}}$ of a proper $C$--set $\mc{S}$ in $\R^n$ that verifies the Minkowski--Lyapunov equation
\begin{equation*}
\forall x\in\R^n,\quad \gauge{\mathcal{S}}{Ax}+\gauge{\mathcal{Q}}{x}= \gauge{\mathcal{S}}{x}.
\end{equation*}
\textbf{Explicit Representation of $D$--sets.}
An explicit representation of a $D$--set $\mc{S}$ is given by the explicit representation of its Minkowski function $\gaugebb{\mc{S}}$. In particular, 
\begin{align*}
x\in\mc{S}\text{ if and only if }&\gauge{\mc{S}}{x}\le 1\text{ so that}\\
\mc{S}=\{x\in\R^n :\ &\gauge{\mc{S}}{x}\le 1\}.
\end{align*}
\textbf{Implicit Representation of $D$--sets.}
An implicit representation of a $D$--set $\mc{S}$ is given by an implicit representation of its Minkowski function $\gaugebb{\mc{S}}$. \emph{The implicit representations of $\mc{S}$ and $\gaugebb{\mc{S}}$ do not require $\mc{S}$ and $\gaugebb{\mc{S}}$ to be explicitly computed}, as exemplified by a relatively direct variation of~\cite[Ch.~7, Sec.~5, Theorem~5]{berge:1963}. 
\begin{theorem}
\label{thm:01.01}
Let $\{\mc{S}_i\ :\ i\in\mc{I}\}$ be a finite collection  of $D$--sets in $\R^n$. The set 
\begin{equation*}
\mc{S}=\bigcap_{i\in\mc{I}}\mc{S}_i
\end{equation*}
is a $D$--set in $\R^n$. Furthermore,
\begin{align*}
\forall x\in\R^n,\quad \gauge{\mc{S}}{x}=&\max_{i\in\mc{I}}\gauge{\mc{S}_i}{x}\text{ and}\\
x\in\mc{S}\text{ if and only if }&\max_{i\in\mc{I}}\gauge{\mc{S}_i}{x}\le 1.
\end{align*}
\end{theorem}

The implicit representations of the Minkowski function $\gaugebb{\mc{S}}$ and its generator set $\mc{S}=\bigcap_{i\in\mc{I}}\mc{S}_i$ are given by
\begin{align*}
x\mapsto &\max_{i\in\mc{I}}\gauge{\mc{S}_i}{x}\text{ and}\\
\mc{S}=\{x\in\R^n :\ &\max_{i\in\mc{I}}\gauge{\mc{S}_i}{x}\le 1\}.
\end{align*}
\textbf{Evaluation of Minkowski Functions.} Any proper $C$--polytopic set $\mc{P}$ has an irreducible representation (in which $\mc{I}_\mc{P}$ is a finite index set and each $p_i\in\R^n$)
\begin{align*}
\mc{P}=\{x\in \R^n\ :\ \forall i\in\mc{I}_\mc{P},\ p_i^Tx\le 1\}.
\end{align*}
Any proper $C$--ellipsoidal set $\mc{E}$ centered at the origin has a representation (in which $E\in\Rnm{n}{n}$ with $E=E^T\succ 0$)
\begin{equation*}
\mc{E}=\{x\in \R^n\ :\ \sqrt{x^TEx}\le 1\}.
\end{equation*}
The evaluation of $\gaugebb{\mc{P}}$ or $\gaugebb{\mc{E}}$ is highly efficient, as
\begin{equation*}
\forall x\in\R^n,\ \gauge{\mc{P}}{x}=\max_{i\in\mc{I}_\mc{P}} p_i^Tx\text{ and }\gauge{\mc{E}}{x}=\sqrt{x^TEx}.
\end{equation*}
The Minkowski function of the intersection of finitely many proper $C$-- polytopic and/or ellipsoidal sets can be also evaluated efficiently, since, as stated in Theorem~\ref{thm:01.01}, 
\begin{equation*}
\forall x\in\R^n,\ \gauge{\bigcap_{i\in\mc{I}}\mc{S}_i}{x}=\max_{i\in\mc{I}}\gauge{\mc{S}_i}{x}.
\end{equation*}
\textbf{Proofs.} The proofs of all formal statements made in this paper are provided in the Appendix.

\section{Alternative Characterization}
\label{sec:02}
The map $\mc{G}\bb$ defined, for subsets  $\mc{S}$ of $\R^n$, by
\begin{align*}
\mc{G}(\mc{S}):=\{x\in\R^n\ :\ &\exists\gamma\in[0,1]\text{ such that }\nonumber\\
&Ax\in \gamma\mc{S}\text{ and }x\in (1-\gamma)\mc{Q}\},
\end{align*}
 and its post fixed points (i.e., sets such that $\mc{S}\subseteq \mc{G}(\mc{S})$) play a crucial role in deriving a novel, alternative, characterization of Minkowski--Lyapunov functions.  
\begin{theorem}
\label{thm:02.01}
Let $A\in \Rnm{n}{n}$ and let $\mc{Q}$ be a proper $C$--set in $\R^n$. $(i)$ A proper $C$--set $\mc{S}$ in $\R^n$ is such that
\begin{equation*}
\forall x\in\R^n,\quad \gauge{\mathcal{S}}{Ax}+\gauge{\mathcal{Q}}{x}\le \gauge{\mathcal{S}}{x}
\end{equation*}
if and only if 
\begin{equation*}
\mc{S}\subseteq \mc{G}(\mc{S}).
\end{equation*}
$(ii)$ A proper $C$--set $\mc{S}$ in $\R^n$ is such that $\mc{S} \subseteq \mc{G}(\mc{S})$ if there exists a scalar $\gamma\in[0,1]$ such that
\begin{equation*}
A\mc{S}\subseteq \gamma \mc{S}\text{ and }\mc{S}\subseteq (1-\gamma)\mc{Q}.
\end{equation*} 
\end{theorem}
Likewise, the map $\mc{G}\bb$ and its maximal fixed point (in the sense that $\mc{S}=\mc{G}(\mc{S})$) are of paramount importance for obtaining a novel, alternative, characterization of the fundamental Minkowski--Lyapunov function.
\begin{theorem}
\label{thm:02.02}
Let $A\in \Rnm{n}{n}$ and let $\mc{Q}$ be a proper $C$--set in $\R^n$. $(i)$ A proper $C$--set $\mc{S}$ in $\R^n$ is such that
\begin{equation*}
\forall x\in\R^n,\quad \gauge{\mathcal{S}}{Ax}+\gauge{\mathcal{Q}}{x}= \gauge{\mathcal{S}}{x}
\end{equation*}
if and only if $\mc{S}$ is the maximal set with respect to set inclusion such that
\begin{equation*}
\mc{S}= \mc{G}(\mc{S}).
\end{equation*} 
$(ii)$ The limit, with respect to the Hausdorff distance, say $\mc{S}$, of the set sequence $\{\mc{S}_k\}_{k\ge 0}$, generated, for all integers $k\ge 0$, by
\begin{equation*}
\mc{S}_{k+1}=\mc{G}(\mc{S}_k)\text{ with } \mc{S}_0=\mc{Q},
\end{equation*}
 is the maximal set with respect to set inclusion such that $\mc{S}=\mc{G}(\mc{S})$. The limit $\mc{S}$ is a $C$--set in $\R^n$. Furthermore, the limit $\mc{S}$ is a proper $C$--set in $\R^n$ if and only if $\rho(A)<1$. 
\end{theorem}

\section{Implicit Representation and Computation}
\label{sec:03}
Theorem~\ref{thm:03.01} identifies a dynamically compatible parametrization of Minkowski--Lyapunov functions $\gaugebb{\mc{S}}$.
\begin{theorem}
\label{thm:03.01}
Let $A\in \Rnm{n}{n}$ be a strictly stable matrix and let $\mc{Q}$ be a proper $C$--set in $\R^n$. 
 For all $\gamma \in (\rho(A),1)$, there exists a finite integer $k>0$ such that
\begin{equation*}
(\gamma^{-1}A)^k\mc{Q}\subseteq \mc{Q}.
\end{equation*}
Furthermore, for all such scalars $\gamma \in (\rho(A),1)$ and integers $k>0$, the set 
\begin{equation*}
\mc{S}=(1-\gamma)\bigcap_{i=0}^{k-1}\left((\gamma^{-1}A)^{-i}\mc{Q}\right)
\end{equation*}
is a proper $C$--set in $\R^n$ such that $\mc{S} \subseteq \mc{G}(\mc{S})$. 
\end{theorem}
Theorem~\ref{thm:03.02} enables an efficient utilization of these parameterized Minkowski--Lyapunov functions $\gaugebb{\mc{S}}$.  
\begin{theorem}
\label{thm:03.02}
Let $\{M_i\in \Rnm{n}{n}\ :\ i\in\mc{I}\}$  and  $\{\mc{S}_i\ :\ i\in\mc{I}\}$ be finite collections of matrices and proper $C$--sets in $\R^n$, respectively. The set 
\begin{equation*}
\mc{S}=\bigcap_{i\in\mc{I}} M_i^{-1}\mc{S}_i
\end{equation*} 
is a proper $D$--set in $\R^n$, which is a proper $C$--set in $\R^n$ when it is bounded. Furthermore,
\begin{align*}
\forall x\in\R^n,\quad \gauge{\mc{S}}{x}=&\max_{i\in\mc{I}} \gauge{\mc{S}_i}{M_ix}\text{ and}\\
x\in \mc{S}\text{ if and only if }&\max_{i\in\mc{I}} \gauge{\mc{S}_i}{M_ix}\le 1.
\end{align*}
\end{theorem}
\begin{remark}
\label{rem:03.01}
\emph{
Evidently, in light of Theorems~\ref{thm:03.01} and~\ref{thm:03.02}, with $\mc{I}=\{0,1,\ldots,k-1\}$, 
\begin{align*}
\forall x\in\R^n,\ \gauge{\mc{S}}{x}=(1-\gamma)^{-1}&\max_{i\in\mc{I}} \gauge{\mc{Q}}{(\gamma^{-1}A)^ix}\text{ and}\\
x\in\mc{S}\text{ if and only if }(1-\gamma)^{-1}&\max_{i\in\mc{I}} \gauge{\mc{Q}}{(\gamma^{-1}A)^ix}\le 1.
\end{align*}
Hence, the implicit representations of Minkowski--Lyapunov functions $\gaugebb{\mc{S}}$ and their generator sets $\mc{S}$ characterized in Theorem~\ref{thm:03.01} are given, respectively, by
\begin{align*}
x\mapsto (1-\gamma)^{-1}&\max_{i\in\mc{I}} \gauge{\mc{Q}}{(\gamma^{-1}A)^ix}\text{ and}\\
\mc{S}=\{x\in\R^n :\ (1-\gamma)^{-1}&\max_{i\in\mc{I}} \gauge{\mc{Q}}{(\gamma^{-1}A)^ix}\le 1\}.
\end{align*}
}
\end{remark}
We identify one more dynamically compatible parametrization of Minkowski--Lyapunov functions $\gaugebb{\mc{S}}$. First, we recall that the polar set $\mc{X}^*$ of a set $\mc{X}$ in $\R^n$ with $0\in \mc{X}$ is 
\begin{equation*}
\mc{X}^*:= \{y\in\R^n\ :\ \forall x\in\mc{X},\ y^T x\le 1\}.
\end{equation*}
\begin{theorem}
\label{thm:03.03}
Let $A\in \Rnm{n}{n}$ be a strictly stable matrix and let $\mc{Q}$ be a proper $C$--set in $\R^n$. 
For all $\gamma \in (0,1)$, there exists a finite integer $k>0$ such that
\begin{equation*}
(A^T)^k\mc{Q}^*\subseteq \gamma \mc{Q}^*. 
\end{equation*}
Furthermore, for all such scalars $\gamma \in (0,1)$ and integers $k>0$, the set 
\begin{equation*}
\mc{S}=\left((1-\gamma)^{-1}\bigoplus_{i=0}^{k-1}(A^T)^{i}\mc{Q}^*\right)^*
\end{equation*}
is a proper $C$--set in $\R^n$ such that $\mc{S} \subseteq \mc{G}(\mc{S})$.
\end{theorem}
We also enable an efficient use of these parameterized Minkowski--Lyapunov functions $\gaugebb{\mc{S}}$.  
\begin{theorem}
\label{thm:03.04}
Let $\{M_i\in \Rnm{n}{n}\ :\ i\in\mc{I}\}$  and  $\{\mc{S}_i\ :\ i\in\mc{I}\}$ be finite collections of matrices and proper $C$--sets in $\R^n$, respectively. The set
\begin{equation*}
\mc{S}=\left(\bigoplus_{i\in\mc{I}} M_i^T\mc{S}_i^*\right)^*
\end{equation*}
is a proper $D$--set in $\R^n$, which is a proper $C$--set in $\R^n$ when it is bounded. Furthermore, 
\begin{align*}
\forall x\in\R^n,\quad \gauge{\mc{S}}{x}=&\sum_{i\in\mc{I}} \gauge{\mc{S}_i}{M_ix}\text{ and}\\
x\in \mc{S}\text{ if and only if }&\sum_{i\in\mc{I}} \gauge{\mc{S}_i}{M_ix}\le 1.
\end{align*}
\end{theorem}
\begin{remark}
\label{rem:03.02}
\emph{
Clearly, in view of  Theorems~\ref{thm:03.03} and~\ref{thm:03.04}, with $\mc{I}=\{0,1,\ldots,k-1\}$, 
\begin{align*}
\forall x\in\R^n,\ \gauge{\mc{S}}{x}=(1-\gamma)^{-1}&\sum_{i\in\mc{I}} \gauge{\mc{Q}}{A^ix}\text{ and}\\
x\in\mc{S}\text{ if and only if }(1-\gamma)^{-1}&\sum_{i\in\mc{I}} \gauge{\mc{Q}}{A^ix}\le 1.
\end{align*}
Thus, the implicit representations of the Minkowski--Lyapunov functions $\gaugebb{\mc{S}}$ and their generator sets $\mc{S}$ characterized in Theorem~\ref{thm:03.03} are given, respectively, by
\begin{align*}
x\mapsto (1-\gamma)^{-1}&\sum_{i\in\mc{I}} \gauge{\mc{Q}}{A^ix}\text{ and}\\
\mc{S}=\{x\in\R^n :\ (1-\gamma)^{-1}&\sum_{i\in\mc{I}} \gauge{\mc{Q}}{A^ix}\le 1\}.
\end{align*}
}
\end{remark}
\begin{remark}
\label{rem:03.03}
\emph{
The pointwise evaluation of the implicit representations $x\mapsto (1-\gamma)^{-1}\max_{i\in\mc{I}} \gauge{\mc{Q}}{(\gamma^{-1}A)^ix}$ and $x\mapsto (1-\gamma)^{-1}\sum_{i\in\mc{I}} \gauge{\mc{Q}}{A^ix}$ of Minkowski--Lyapunov functions $\gaugebb{\mc{S}}$ characterized in Theorems~\ref{thm:03.01} and~\ref{thm:03.03}, respectively, is efficient in arbitrary finite dimensions for a rich variety of the proper $C$--sets $Q$. It is very simple and highly efficient when $\mc{Q}$ is a proper $C$--set that is either polytopic or ellipsoidal set or the intersection of polytopic and/or ellipsoidal sets. (See remarks on the evaluation of Minkowski functions in Section~\ref{sec:01}.) In particular,  depending on the considered case, for a given $x$, one generates the sequence of points $\{(\gamma^{-1}A)^ix\}_{i=0}^{k-1}$ or $\{A^ix\}_{i=0}^{k-1}$  and evaluates the sequence of values $\{(1-\gamma)^{-1}\gauge{\mc{Q}}{(\gamma^{-1}A)^ix}\}_{i=0}^{k-1}$ or $\{(1-\gamma)^{-1}\gauge{\mc{Q}}{A^ix}\}_{i=0}^{k-1}$, and then computes the maximum or the sum of the latter sequence.
}
\end{remark}
\begin{remark}
\label{rem:03.04}
\emph{
In order to utilize the implicit representations of Minkowski--Lyapunov functions $\gaugebb{\mc{S}}$ characterized in Theorems~\ref{thm:03.01} and~\ref{thm:03.03}, all what is needed is to detect an integer (possibly the minimal integer) $k$ for which the conditions postulated in Theorems~\ref{thm:03.01} and~\ref{thm:03.03} hold true. These conditions take a generic form $M^k\mc{X}\subseteq \mc{X}$ for a strictly stable matrix $M\in\Rnm{n}{n}$ and a proper $C$--set $\mc{X}$ in $\R^n$. Such a set inclusion can be also handled efficiently in arbitrary finite dimensions for a rich variety of the proper $C$--sets $\mc{X}$. In particular, for a proper $C$--ellipsoidal set $\mc{X}=\{x\ :\ \sqrt{x^TXx}\le 1\}$, such a set inclusion is equivalent to $(M^k)^TXM^k-X\preceq 0$, which can be checked, for instance, by evaluating the eigenvalues of the matrix $(M^k)^TXM^k-X$. Likewise, for a proper $C$--polytopic set $\mc{X}$ with an irreducible representation $\mc{X}=\{x\ :\ \forall i\in \mc{I}_{\mc{X}},\ x_i^Tx\le 1\}$, such a set inclusion holds true if and only if, for all $i\in \mc{I}_{\mc{X}}$, $\support{\mc{X}}{(M^k)^Tx_i}\le 1$, where $\supportbb{\mc{X}}$ is the support function  (defined in the next page) and which can be checked, for example, by solving
$\operatorname{cardinality}(\mc{I}_{\mc{X}})$ linear programming problems. The above observations can be combined so as to address (via sufficiency) the case when $\mc{X}$ is the intersection of finitely many proper $C$-- polytopic and/or ellipsoidal sets.
}
\end{remark}
Theorem~\ref{thm:02.02} characterizes the generator set $\mc{S}$ of the fundamental Minkowski--Lyapunov function $\gaugebb{\mc{S}}$ as the maximal fixed point of the map $\mc{G}\bb$, which is the limit, with respect to the Hausdorff distance~\cite{klein:thompson:1984}, of the set sequence $\{\mc{S}_k\}_{k\ge 0}$ generated by the set recursion specified in Theorem~\ref{thm:02.02}$(ii)$. This characterization and the corresponding set recursion can be seen as a constructive utilization of the Tarski fixed point theorem~\cite{tarski:1955} and the Kleene--like iteration~\cite{kleene:1952}. Note that, for any strictly stable matrix $A\in \Rnm{n}{n}$ and any proper $C$--set $\mc{Q}$ in $\R^n$, the sets $\mc{S}_k,\ k\ge 0$ are proper $C$--sets in $\R^n$ such that $\mc{S}_{k+1}\subseteq \mc{S}_k$ and their limit $\mc{S}$ is a proper $C$--set in $\R^n$. The limit $\mc{S}$ is finitely determined when $\mc{S}_k\subseteq \mc{S}_{k+1}$ for an integer $k\ge 0$, in which case $\mc{S}=\mc{S}_k=\mc{S}_{k+1}$.  This set recursion is the iterative evaluation of the map  $\mc{G}\bb$ at $\mc{Q}$. Its worst case computational complexity can be considerable in general. Consequently, its  numerically plausible implementation should take advantage of any available structure of the matrix $A$ and proper $C$--set $\mc{Q}$. 

\section{Intrinsic Duality Implications}
\label{sec:04}

By~\cite[Theorem~1]{rakovic:2020.a}, the Minkowski function $\gaugebb{\mc{S}}$ of a proper $C$--set $\mc{S}$ in $\R^n$ verifies the Minkowski--Lyapunov inequality if and only if the polar set $\mc{Z}=\mc{S}^*$ is such that
\begin{equation*}
A^T\mc{Z}\oplus \mc{W}\subseteq \mc{Z}\text{ with }\mc{W}=\mc{Q}^*.
\end{equation*}
This set inclusion is a necessary and sufficient condition for a set $\mc{Z}$ to be a robust positively invariant set~\cite{kolmanovsky:gilbert:1998} for the polar linear dynamics 
\begin{equation*}
z^+=A^Tz+w\text{ with }w\in \mc{W}=\mc{Q}^*,
\end{equation*}
for which  $z\in \R^n$, $w\in \R^n$ and $z^+\in \R^n$ are the polar current state and disturbance and polar successor state. By~\cite[Theorem~3]{rakovic:2020.a}, the Minkowski function $\gaugebb{\mc{S}}$ of a proper $C$--set $\mc{S}$ in $\R^n$ verifies the Minkowski--Lyapunov equation if and only if the polar set $\mc{Z}=\mc{S}^*$ is such that 
\begin{equation*}
A^T\mc{Z}\oplus \mc{W}= \mc{Z}\text{ with }\mc{W}=\mc{Q}^*.
\end{equation*}
This  fixed point set equation is, within the considered setting, a necessary and sufficient condition~\cite{artstein:rakovic:2008} for a set $\mc{Z}$ to be the minimal (nonempty and compact) robust positively invariant set~\cite{kolmanovsky:gilbert:1998} for the polar linear dynamics. 

The preceding facts and Theorems~\ref{thm:02.01}$(i)$ and~\ref{thm:02.02}$(i)$ yield directly alternative characterizations of robust positively invariant proper $C$--sets and the minimal robust positively invariant set (over the space of nonempty compact subsets of $\R^n$) for the polar linear dynamics.
\begin{corollary}
\label{cor:04.01}
Let $A\in \Rnm{n}{n}$ and let $\mc{Q}$ be a proper $C$--set in $\R^n$. $(i)$ A proper $C$--set $\mc{Z}$ in $\R^n$ is such that
\begin{equation*}
A^T\mc{Z}\oplus \mc{W} \subseteq \mc{Z}\text{ with }\mc{W}=\mc{Q}^*
\end{equation*}
if and only if its polar set $\mc{Z}^*$ is such that $\mc{Z}^*\subseteq \mc{G}(\mc{Z^*})$.
$(ii)$ A proper $C$--set $\mc{Z}$ in $\R^n$ is such that
\begin{equation*}
A^T\mc{Z}\oplus \mc{W} = \mc{Z}\text{ with }\mc{W}=\mc{Q}^*
\end{equation*} 
if and only if  its polar set $\mc{Z}^*$ is the maximal set with respect to set inclusion such that $\mc{Z}^*= \mc{G}(\mc{Z^*})$.
\end{corollary}

Hence, the polar sets $\mc{Z}=\mc{S}^*$ of the proper $C$--sets $\mc{S}$ characterized in Theorems~\ref{thm:03.01} and~\ref{thm:03.03} are robust positively invariant proper $C$--sets for the polar linear dynamics. We recall that the support function of a nonempty closed convex set $\mc{X}$ in $\R^n$ is given, for all $y\in\R^n$, by
\begin{equation*}
\support{\mc{X}}{y}:= \sup_x \{y^Tx\ :\ x\in \mc{X}\}.
\end{equation*} 
By the virtue of~\cite[Theorems~1.6.1 and~1.7.6]{schneider:1993}, the support functions $\supportbb{\mc{Z}}$ of the  polar sets $\mc{Z}=\mc{S}^*$ of the proper $C$--sets $\mc{S}$ characterized in Theorems~\ref{thm:03.01} and~\ref{thm:03.03} satisfy $\support{\mc{Z}}{x}=\gauge{\mc{Z}^*}{x}=\gauge{\mc{S}}{x}$ for all $x\in\R^n$, and, thus, are given, respectively, by 
\begin{align*}
\forall x\in\R^n,\quad \support{\mc{Z}}{x}=(1-\gamma)^{-1}&\max_{i\in\mc{I}} \gauge{\mc{Q}}{(\gamma^{-1}A)^ix}\text{ and}\\
\forall x\in\R^n,\quad \support{\mc{Z}}{x}=(1-\gamma)^{-1}&\sum_{i\in\mc{I}} \gauge{\mc{Q}}{A^ix},
\end{align*}
which reveals directly their efficient implicit representations, and which when substituted in
\begin{equation*}
\mc{Z}=\{x\in\R^n :\ \forall y\in\R^n,\ y^Tx\le \support{\mc{Z}}{y}\}
\end{equation*}
yields the efficient implicit representations of the related robust positively invariant proper $C$--sets $\mc{Z}=\mc{S}^*$.

\section{Closing Remarks and Numerical Experience}
\label{sec:05}

Even the computation of polyhedral positively invariant sets~\cite{bitsoris:1988} and Lyapunov functions are relevant and active research topics. For a plethora of theoretical and computational contributions, see, for instance,~\cite{blanchini:1999,blanchini:miani:2008,giesl:hafstein:2015} and references therein. In relation to the existing methods, and as discussed more formally in Remarks~\ref{rem:03.03} and~\ref{rem:03.04}, the implicit representations of  Minkowski--Lyapunov functions identified in Theorems~\ref{thm:03.01} and~\ref{thm:03.03} (and, by intrinsic duality, the corresponding robust positively invariant sets) are numerically potent in arbitrary finite dimensions. Furthermore, their structure is flexible since it is neither restricted to polytopic nor ellipsoidal proper $C$--sets. 
  
  Table~1 summarizes the outcome of a numerical test with a sample of randomly generated strictly stable matrices 
\begin{table}[!hbt]    
\begin{center}
\resizebox{.4875\textwidth}{!}{%
  \begin{tabular}{c|ccccccc}
  \hline
     $n$ &  $2$ & $3$ & $5$& $8$ & $13$ & $21$ & $34$\\ 
    \hline
   	$\rho(A)$ &   $0.977$ & $0.986$ & $0.982$ & $0.985$ & $0.993$ & $0.978$ & $0.981$\\ 
     $k$ &  $16$ & $19$ & $44$ & $57$ & $100$ & $50$ & $84$\\ 
     $t$ &  $0.1$ & $0.3$ & $0.1$ & $0.2$ & $0.6$ & $0.7$ & $3.5$\\
      \hline
     $n$ &  $55$ & $89$ & $144$& $233$ & $377$ & $610$ & $987$\\ 
    \hline    
   	$\rho(A)$ &   $0.999$ & $0.992$ & $0.999$ & $0.998$ & $0.998$ & $0.995$ & $0.991$\\ 
     $k$ &  $955$ & $270$ & $1272$ & $801$ & $811$ & $624$ & $351$\\ 
     $t$ &  $116$ & $82.9$ & $1430$ & $3308$ & $17275$ & $56512$ & $122139$\\
    \hline
  \end{tabular}
}
Table 1. MATLAB computations for sets $\mc{S}$ of Theorem~\ref{thm:03.01}.
\end{center}
\end{table}
$A\in\Rnm{n}{n}$ for which $\rho(A)\in [0.975,0.999]$, and with $\mc{Q}=\mc{B}^n_\infty:=\{x\in\R^n\ :\ \|x\|_\infty \le 1\}$ and $\gamma=(\rho(A)+1)/2$.
 Table~1 reports the spectral radius $\rho(A)$, the minimal integer $k$ needed to construct the implicit representations of the related  Minkowski--Lyapunov functions $\gaugebb{\mc{S}}$, and the time in milliseconds (ms) needed to compute such an integer $k$ by means of a direct incremental search. The computational times vary from $0.1\ ms$ for $2$--dimensional and $5$--dimensional examples to $122139\ ms$ (i.e., about $2$ minutes and $2$ seconds) for $987$--dimensional example. Clearly, the data reported in Table~1  furnishes strong evidence of the asserted numerical potency of Minkowski--Lyapunov functions characterized in Theorem~\ref{thm:03.01}.

Theorem~\ref{thm:02.02} delivers an alternative approach for the computation of the fundamental Minkowski--Lyapunov function and, by intrinsic duality, of the related minimal robust positively invariant set. This approach is also novel and, more importantly, it has a potential to be numerically plausible within more structured settings. 
\begin{figure}[!h]
\includegraphics[width=0.4875\textwidth]{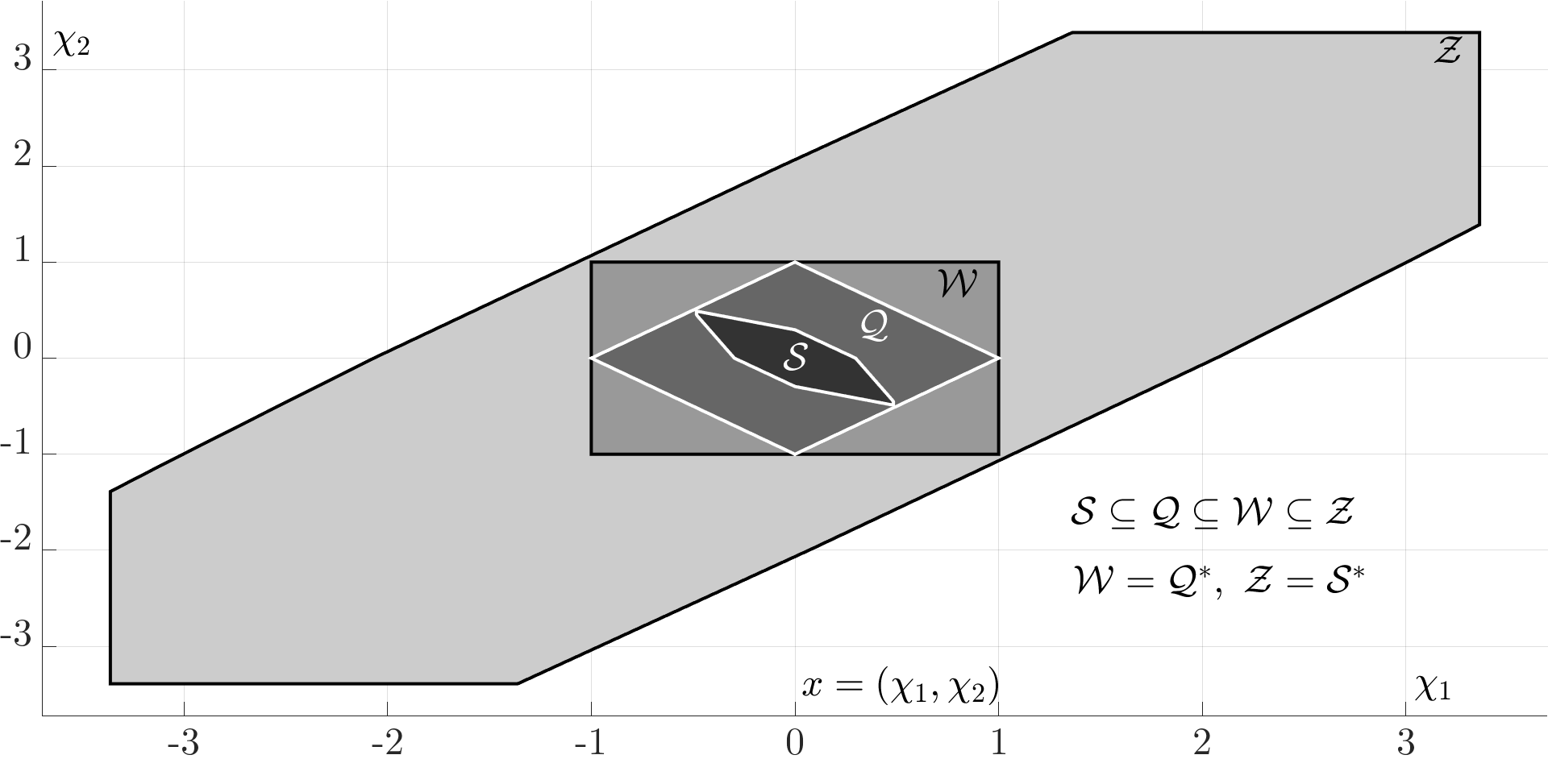}  
\caption{The Sets $\mc{S}$, $\mc{Q}$, $\mc{W}=\mc{Q}^*$ and $\mc{Z}=\mc{S}^*$.}
\label{fig:04.01}
\end{figure}
Figure~\ref{fig:04.01} illustrates the related polyhedral computations for an example, in which the entries of the matrix $A$ are $a_{(1,1)}=1$, $a_{(1,2)}=1$, $a_{(2,1)}=-0.72$, and $a_{(2,2)}=-0.7$ so that $\rho(A)=0.2$, and $\mc{Q}=\mc{B}_1:=\{x\ :\ \|x\|_1\le 1\}$ so that $\mc{W}=\mc{Q}^*=\mc{B}_\infty:=\{x\ :\ \|x\|_\infty \le 1\}$.  
The set iteration proposed in Theorem~\ref{thm:02.02} generates its (numerical) limit $\mc{S}$ in $k=10$ steps, which is a proper $C$--polytopic set and the maximal set such that $\mc{S}=\mc{G}(\mc{S})$. By Corollary~\ref{cor:04.01}, the corresponding polar set $\mc{Z}=\mc{S}^*$ is the minimal proper $C$--polytopic set such that $A^T\mc{Z}\oplus \mc{W}=\mc{Z}$. 

In terms of future research, a more dedicated study, with a primary focus on the related computational aspects, of a setting with a more structured matrix $A$ and the proper $C$--set $\mc{Q}$, would further enhance a  practical utilization of Minkowski--Lyapunov functions and the fundamental Minkowski--Lyapunov function. Likewise, the study of the efficient implicit representations of the fundamental Minkowski--Lyapunov function (and, by inherent duality, the corresponding minimal robust positively invariant set) would be also of much interest. 

\textbf{Acknowledgement.} This is an extended version of the accepted Automatica technical communiqu\'{e} $20$--$1623$, which was initially submitted on November $25$, $2020$, and which is to be published by Automatica. The author is grateful to the Editor, Associate Editor and Referees for timely reviews and very constructive comments.

\section*{Appendix: Proofs}

\textbf{A. Proof of Theorem~\ref{thm:01.01}.} 
By definition, the set $\mc{S}=\bigcap_{i\in\mc{I}} \mc{S}_i$ is at least a $D$--set in $\R^n$. By~\cite[Ch.~7, Sec.~5, Theorem~5]{berge:1963}, for all $x\in \R^n$, $\gauge{\mc{S}}{x}=\max_{i\in\mc{I}}\gauge{\mc{S}_i}{x}$. Hence, the claim.

\textbf{B. Proof of Theorem~\ref{thm:02.01}.}
$(i)$ Since $\mc{S}$ and $\mc{Q}$ are proper $C$--sets and the Minkowski function $\gaugebb{\mc{X}}$ of a proper $C$--set $\mc{X}$ is positively homogeneous of the first degree, it suffices to establish the claim for an arbitrary $x$ such that $\gauge{\mc{S}}{x}=1$, so in this proof we consider such an $x$.

 First, let $\mc{S}$ be such that $\gauge{\mc{S}}{Ax}+\gauge{\mc{Q}}{x}\le \gauge{\mc{S}}{x}$. Let $\gamma:=\gauge{\mc{S}}{Ax}$ so that $\gamma \in [0,1]$ and $Ax\in \gamma \mc{S}$. Also, $\gauge{\mc{Q}}{x}\le \gauge{\mc{S}}{x}-\gauge{\mc{S}}{Ax}=1-\gamma$ so that $x\in (1-\gamma)\mc{Q}$. Since  $\gamma \in [0,1]$, $Ax\in \gamma \mc{S}$ and $x\in (1-\gamma)\mc{Q}$, it follows that $x\in \mc{G}(\mc{S})$. Hence, $\mc{S}\subseteq \mc{G}(\mc{S})$. 
 
 Second, let $\mc{S}$ be such that $\mc{S}\subseteq \mc{G}(\mc{S})$. Then there exists a $\gamma \in [0,1]$ such that $Ax\in \gamma \mc{S}$ and $x\in (1-\gamma)\mc{Q}$. In turn, $\gauge{\mc{S}}{Ax}\le \gamma$ and $\gauge{\mc{Q}}{x}\le 1-\gamma$. Consequently, $\gauge{\mc{S}}{Ax}+\gauge{\mc{Q}}{x}\le \gamma+(1-\gamma)=1=\gauge{\mc{S}}{x}$. Hence, $\mc{S}$ is such that $\gauge{\mc{S}}{Ax}+\gauge{\mc{Q}}{x}\le \gauge{\mc{S}}{x}$. 
 
 $(ii)$ Since $\mc{S}$ and $\mc{Q}$ are proper $C$--sets, the claimed fact follows from~$(i)$.

\textbf{C. Proof of Theorem~\ref{thm:02.02}.}
$(i)$ This fact follows from Theorem~\ref{thm:02.01}$(i)$,~\cite[Theorems~3]{rakovic:2020.a} and the postulated maximality of the proper $C$--set $\mc{S}$. 

$(ii)$ The collection of all convex subsets of $\R^n$ is a complete lattice under the natural partial ordering corresponding to the set inclusion~\cite{rockafellar:1970}, and, by its definition, $\mc{G}\bb$ maps convex subsets of $\R^n$ into  convex subsets of $\R^n$ and it is monotone (i.e. $\mc{X}\subseteq \mc{Y}$ implies that $\mc{G}(\mc{X})\subseteq \mc{G}(\mc{Y})$). Thus, all postulates of the Tarski fixed point theorem~\cite{tarski:1955} are satisfied, and the Tarski fixed point theorem  guarantees the existence of  the maximal fixed point $\mc{S}$ of the map $\mc{G}\bb$ over the collection of convex subsets of $\R^n$. Since, for any subset $\mc{X}$ of $\R^n$, $\mc{G}(\mc{X})\subseteq \mc{Q}$ and $\mc{G}\bb$ maps (proper) $C$--sets in $\R^n$ into (proper) $C$--sets in $\R^n$, the maximal fixed point $\mc{S}$ of $\mc{G}\bb$ is guaranteed to be a $C$--set in $\R^n$, which is a proper $C$--set in $\R^n$ if and only if $\rho(A)<1$. Namely,  when the maximal fixed point $\mc{S}$ of $\mc{G}\bb$ is a proper $C$--set in $\R^n$, Minkowski--Lyapunov function $\gaugebb{\mc{S}}$ verifies the strict stability of the matrix $A$. Likewise, when the matrix $A$ is strictly stable there exists a Minkowski--Lyapunov function $\gaugebb{\mc{R}}$ generated by a proper $C$--set $\mc{R}$ in $\R^n$ so that $\mc{R}\subseteq\mc{G}(\mc{R})\subseteq \mc{S}=\mc{G}(\mc{S})$ and the maximal fixed point $\mc{S}$ of $\mc{G}\bb$  is a proper $C$--set in $\R^n$. In either case, $\mc{S}\subseteq \mc{Q}$ and the maximal fixed point $\mc{S}$ of $\mc{G}\bb$ is the limit, with respect to the Hausdorff distance~\cite{klein:thompson:1984}, of the set sequence $\{\mc{S}_k\}_{k\ge 0}$ generated by the considered set recursion. 

\textbf{D. Proof of Theorem~\ref{thm:03.01}.}
Since $\rho(A)\in [0,1)$ and $\gamma \in (\rho(A),1)$, $\rho(\gamma^{-1}A)\in [0,1)$. Hence, since $\rho(\gamma^{-1}A)\in [0,1)$ and  $\mc{Q}$ is a proper $C$--set in $\R^n$, there exists a finite integer $k>0$ such that $(\gamma^{-1}A)^k\mc{Q}\subseteq \mc{Q}$.  For any such $k$ and $\gamma$, by definition, $\mc{S}\subseteq (1-\gamma) \mc{Q}$ is a proper $C$--set in $\R^n$. Since $(\gamma^{-1}A)^k\mc{Q}\subseteq \mc{Q}$, $\mc{Q}\subseteq (\gamma^{-1}A)^{-k}\mc{Q}$ and, thus,
\begin{equation*}
\mc{S}=(1-\gamma)\bigcap_{i=0}^{k-1}(\gamma^{-1}A)^{-i}\mc{Q}\subseteq (1-\gamma)\bigcap_{i=1}^{k}(\gamma^{-1}A)^{-i}\mc{Q},
\end{equation*}
and the proof is concluded by noting that $A\mc{S}\subseteq \gamma \mc{S}$, as
\begin{align*}
A\mc{S}
&\subseteq \gamma (\gamma^{-1}A)(1-\gamma)\bigcap_{i=1}^{k}(\gamma^{-1}A)^{-i}\mc{Q}\\
&\subseteq \gamma (1-\gamma)\bigcap_{i=0}^{k-1}(\gamma^{-1}A)^{-i}\mc{Q}=\gamma \mc{S}.
\end{align*}
\textbf{E. Proof of Theorem~\ref{thm:03.02}.}
By definition, the set $\mc{S}=\bigcap_{i\in\mc{I}} M_i^{-1}\mc{S}_i$ is at least a proper $D$--set in $\R^n$ and a proper $C$--set in $\R^n$ when it is bounded. By~\cite[Ch.~7, Sec.~5, Theorem~5]{berge:1963}, for all $x\in \R^n$, $\gauge{\mc{S}}{x}=\max_{i\in\mc{I}}\gauge{M_i^{-1}\mc{S}_i}{x}$. By~\cite[Corollary~16.3.2]{rockafellar:1970}, for all $x\in \R^n$ and all $i\in\mc{I}$, $\gauge{M_i^{-1}\mc{S}_i}{x}=\gauge{\mc{S}_i}{M_ix}$. Hence,
\begin{align*}
\forall x\in\R^n,\quad \gauge{\mc{S}}{x}=&\max_{i\in\mc{I}} \gauge{\mc{S}_i}{M_ix}\text{ and}\\
x\in \mc{S}\text{ if and only if }&\max_{i\in\mc{I}} \gauge{\mc{S}_i}{M_ix}\le 1.
\end{align*}
\textbf{F. Proof of Theorem~\ref{thm:03.03}.}
The claimed result follows directly from~\cite[Theorem~1]{rakovic:kerrigan:kouramas:mayne:2004b} and ~\cite[Theorem~1]{rakovic:2020.a}.

\textbf{G. Proof of Theorem~\ref{thm:03.04}.}
By definition, the set $\mc{S}=\left(\bigoplus_{i\in\mc{I}} M_i^T\mc{S}_i^*\right)^*$  is at least a proper $D$--set in $\R^n$ and a proper $C$--set in $\R^n$ when it is bounded. By~\cite[Theorem~14.5]{rockafellar:1970} and the algebra~\cite{rockafellar:1970,schneider:1993} of support functions, for all $x\in \R^n$, 
\begin{align*}
\gauge{\mc{S}}{x}&=\support{\mc{S}^*}{x}=\support{\bigoplus_{i\in\mc{I}} M_i^T\mc{S}_i^*}{x}=\sum_{i\in\mc{I}}\support{M_i^T\mc{S}_i^*}{x}\\
&=\sum_{i\in\mc{I}}\support{\mc{S}_i^*}{M_ix}=\sum_{i\in\mc{I}}\gauge{\mc{S}_i}{M_ix}. 
\end{align*}
Hence,
\begin{align*}
\forall x\in\R^n,\quad \gauge{\mc{S}}{x}=&\sum_{i\in\mc{I}} \gauge{\mc{S}_i}{M_ix}\text{ and}\\
x\in \mc{S}\text{ if and only if }&\sum_{i\in\mc{I}} \gauge{\mc{S}_i}{M_ix}\le 1.
\end{align*} 
\textbf{H. Proof of Corollary~\ref{cor:04.01}.}
The stated facts follow from Theorems~\ref{thm:02.01}$(i)$ and~\ref{thm:02.02}$(i)$ and~\cite[Theorem~1 and~3]{rakovic:2020.a}.

\bibliographystyle{unsrt}
\bibliography{MLFACIC.arXiv}
\end{document}